\documentclass[a4paper,13 pt]{article}
\usepackage [dvips]{graphicx,color}
\normalsize

\begin{document}
\begin{center}
\textbf{PLURIPOLARITY OF GRAPHS OF DENJOY QUASIANALYTIC FUNCTIONS
OF SEVERAL VARIABLES}
\end{center}

\begin{center}
S.A. Imomkulov, Z.Sh. Ibragimov
\end{center}

\begin {abstract}

In this paper we prove pluripolarity of graphs of Denjoy
quasianalytic functions of several variables on the spanning set

\[
T^n = \{z \in {\rm C}^n:\,\,\,\left| {z_1 } \right| = \left| {z_2
} \right| = ... = \left| {z_n } \right| = 1\}.
\]

\end{abstract}

Quasianalytic functions have become the subject of many
investigations in recent years due to their applications in
pluripotential theory and multidimensional complex analysis. The graphs of quasianalytic functions are
closely related to pluripolar sets, which are the main objects of
pluripotential theory. Such relations were studied by several
authors (e.g. J.E. Fornass, K. Diederich [1,2], N. Shcherbina [6],
T. Edlund, B. Joricke [4], D. Coman, N. Levenberg, E. Poletsky
[5], A. Edigarian, J. Wiegerinck [3] and others).

It is easy to show that if $f \in A[a,b]$, then its graph

\[
\Gamma_f = \{(x,f(x)) \in {\rm C}^2:\,x \in [a,b]\}
\]
is pluripolar set in ${\rm C}^2$. In [1] J.E. Fornass and K.
Diederich constructed an example of a ${\rm C}^\infty $ function
$f$ with nonpluripolar graph in ${\rm C}^2$. Recently, D.
Coman, N. Levenberg and E. Poletsky [5] have proved that the graph of
Denjoy quasianalytic functions is pluripolar in ${\rm C}^2$, i.e.
if $f:T \to {\rm C},\,\,\,T = \{\left| z \right| = 1\}$ Denjoy
quasianalytic function, then its graph pluripolar in ${\rm C}^2$.

In this work we consider Denjoy quasianalytic functions of several variables
on the spanning set

\[
T^n = \{z \in {\rm C}^n:\,\,\,\left| {z_1 } \right| = \left| {z_2 } \right|
= ... = \left| {z_n } \right| = 1\}.
\]

Let $M_j $ be a sequence of positive numbers. We denote by $C_{M_j
} (T^n)$ the class of infinity differentiable functions $f \in
C^\infty (T^n)$ satisfying the condition

\[êô
M_j (f) \le \,R^jM_j ,
\]
where $R$ depends on $f$. The class $C_{M_j } (T^n)$ is called the class of Denjoy quasianalytic functions
if for any two functions $f,g \in C_{M_j } (T^n)$ the condition

\[
f^{(\alpha )}(z^0) = g^{(\alpha )}(z^0),\,\,\,\alpha \in {\rm Z}_ + ^n ,
\]

\noindent
at some point $z^0 \in T^n$ implies that

\[
f(z) \equiv g(z).
\]

\textbf{Definition.} \textit{A function }$f \in C^\infty (T^n)$\textit{ is called Denjoy quasianalytic if the class }$C_{M_j (f)} (T^n)$\textit{ is the class of Denjoy quasianalytic functions. }

According to Carleman's theorem [7], a function $f \in C^\infty (T^n)$ is
Denjoy quasianalytic if and only if

\[
\int\limits_1^\infty {\frac{\ln \tau _f (r)}{r^2}} dr = - \infty ,
\]

\noindent
where

\[
\tau _f (r) = \mathop {\inf }\limits_{j \ge 0} \frac{M_j (f)}{r^j}.
\]

Let $f:T^n \to {\rm C}$ be infinity differentiable function with multiple
Fourier series:

\begin{equation}
\label{eq1}
f(e^{i\theta }) = \sum\limits_{k \in {\rm Z}^n} {c_k e^{ik\theta }} ,
\end{equation}

\noindent
where

\[
k = (k_1 ,\,k_2 ,....,k_n ),\,\,\,\,k_j \in {\rm Z},\,\,\,\,j =
1,2,...,n, \quad \theta = (\theta _1 ,\theta _2 ,....,\theta _n
),\,\,\,\,\,\,0 \le \theta _j \le 2\pi ,\,\,\,\,\,\]
\[j =1,2,...,n,   \,\,\,\,\   e^{ik\theta } = e^{ik_1 \theta _1 }e^{ik_2 \theta _2 }
\cdot \cdot \cdot e^{ik_n \theta _n }.
\]

We consider $L_2 $ - norm of partial derivatives:

\begin{equation}
\label{eq2}
\frac{1}{(2\pi )^n}\int\limits_0^{2\pi } {\int\limits_0^{2\pi } { \cdot
\cdot \cdot \int\limits_0^{2\pi } {\left| {\tilde {f}^{(\alpha )}(\theta )}
\right|^2d\theta = \sum\limits_{k \in {\rm Z}^n\backslash \{k:\,k_p =
0,\,\,\alpha _p \ne 0\}} {k^{2\alpha }\left| {c_k } \right|^2} } } } \le
M_j^2 (f),
\quad
\left| \alpha \right| = j,
\end{equation}

\noindent
where $\alpha = (\alpha _1 ,\alpha _2 ,....,\alpha _n ) \in {\rm Z}_ + ^n $,
$k^{2\alpha } = k_1^{2\alpha _1 } k_2^{2\alpha _2 } ....k_n^{2\alpha _n } $;
if $k_p = \alpha _p = 0$, we let $k_p^{2\alpha _p } = 1,\,\,\,\,1 \le p \le
n$.

\textbf{Theorem. }\textit{Let }$f:T^n \to {\rm C}$\textit{ be a Denjoy quasianalytic function. Then its graph }$\Gamma _f $\textit{ pluripolar in }${\rm C}^{n + 1}.$

First, we will prove following proposition.

\textbf{Proposition. }\textit{Let }$f(z) \in C^\infty (T^n)$\textit{ be a function such that }

\begin{equation}
\label{eq3}
\overline {\mathop {\lim }\limits_{m \to \infty } } m^{\textstyle{1 \over {n
+ 1}}}\ln t_m = \infty ,
\end{equation}

\textit{where}

\begin{equation}
\label{eq4}
\ln t_m = \inf \{ - \frac{\ln r^3\tau _f (r)}{nr}:\,\,1 \le r \le m\}.
\end{equation}

\textit{Then its graph }$\Gamma _f $\textit{ - pluripolar.}

\textbf{Proof. }Let $f(z) \in C^\infty (T^n)$ be infinity differentiable
function with Fourier expansion

\[
f(e^{i\theta }) = \sum\limits_{k \in {\rm Z}^n} {c_k e^{ik\theta }}.
\]

We consider the following rational function corresponding to Fourier expansion of $f$:

\begin{equation}
\label{eq5}
L_m (f;z) = \sum\limits_{r = 0}^{m - 1} {\sum\limits_{\beta \in \Lambda }
{a_{m,r} z_1^{\beta _1 r} z_2^{\beta _2 r} ...z_n^{\beta _n r} } } ,
\quad
m = 1,2,....
\end{equation}

\noindent
where

\begin{equation}
\label{eq6}
a_{m,r} = \sum\limits_{l \in {\rm Z}_ + ^n } {\sum\limits_{\beta \in \Lambda
} {c_{\beta r + m\beta l} } } ,
\end{equation}

\[
\begin{array}{l}
 \,\,\,\,\,\Lambda = \left\{ {\beta = (\beta _1 ,\beta _2 ,...,\beta _n
):\,\,\beta _j = \pm 1,\,j = 1,2,...,n} \right\}, \\
 \,\,\,\beta r = (\beta _1 r,\beta _2 r,...,\beta _n r)\,,\,\,\,\,\,m\beta l
= (m\beta _1 l_1 ,m\beta _2 l_2 ,...,m\beta _n l_n ).\,\,\,\, \\
 \end{array}
\]

According to (\ref{eq2}) we have the following estimates for the Foureir
coefficients:

\begin{equation}
\label{eq7}
\left| {c_k } \right| \le \frac{M_j (f)}{\left| {k_1 } \right|^{\alpha _1
}\left| {k_2 } \right|^{\alpha _2 } \cdot \cdot \cdot \left| {k_n }
\right|^{\alpha _n }},\,\,\,\,\left| k \right| \ge 1,\,\,\,\alpha
:\,\,\,\left| \alpha \right| = j \ge 2n.
\end{equation}

$\alpha _p = 0$ if $k_p = 0$ and $\alpha _p \ge 2$ if $k_p \ne 0,\,\,\,\,\,p =
1,2,....,n$; (if $k_p = \alpha _p = 0$, we let $\left| {k_p }
\right|^{\alpha _p } = 1)$.

According to (\ref{eq7}) the series$\sum\limits_{k \in {\rm Z}^n} {\left| {c_k }
\right|} $ is convergent and, consequently, the series (\ref{eq6}) converges
absolutely.

Now we will construct an interpolation trigonometric polynomial $\{L_m
(f;z^0;z)\}$, $m \in {\rm Z}_ + $, where $z^0 = e^{i\theta _0 }$ some fixed
point from $T^n$ and $z = e^{i\theta }$ is an arbitrary point from $T^n$:

\[
L_m (f;z^0;z) = L_m (f;z) + \frac{z_1^m + z_2^m + ... + z_n^m - n}{(z_1^0
)^m + (z_2^0 )^m + ... + (z_n^0 )^m - n}(f(z^0) - L_m (f,z^0)).
\]

Here

\[
L_m (f;z^0;z^{(l)}) = f(z^{(l)}),\,\,\,
\]

\noindent
where

\[
z^{(l)} = (e^{\textstyle{{2\pi il_1 } \over m}},e^{\textstyle{{2\pi il_2 }
\over m}},...,e^{\textstyle{{2\pi il_n } \over m}}),
\]

$l \in J_m \, = \,\left\{ {l = (l_1 ,l_2 ,...,l_n ) \in {\rm Z}_ +
^n ,\,\,1 \le l_p \le m} \right\}$. The number of such points on
$T^n$ are equal $m^n$. In case $(z_1^0 )^m + (z_2^0 )^m + ... +
(z_n^0 )^m = n$, the value of $L_m (f;z^0;z)$ is defined as $L_m
(f;z^0;z) = L_m (f;z)$. The function
$L_m (f;z^0;z)$ interpolates $f$ at the points $z^0,z^{(l)},\,l
\in J_m $.

Now we show that there exists a constant $C_f $, which depends on $f$, such
that for each $m \ge 1$, $z^0 \in T^n$ and $t > 1$ we have

\begin{equation}
\label{eq8}
\left| {L_m (f;z^0;z)} \right| \le C_f (1 + \sum\limits_{r = 1}^m {r\tau _f
(r)t^{nr})}
\end{equation}

\noindent
for all $z$: $\frac{1}{t} \le \left| {z_p } \right| \le t,\,\,\,\,p =
1,2,...,n$.
In fact, according to inequality (\ref{eq7}) we have
$$
\left| {a_{m,r} } \right| \le \sum\limits_{l \in {\rm Z}_ + ^n }
{\sum\limits_{\beta \in \Lambda } {\left| {c_{\beta r + m\beta l}
} \right| \le 2^nM_j (f)\sum\limits_{l \in {\rm Z}_ + ^n }
{\frac{1}{\left| {r + ml_1 } \right|^{\alpha _1 }\left| {r + ml_2
} \right|^{\alpha _2 }...\left| {r + ml_n } \right|^{\alpha _n }}}
} } \le \frac{4^nM_j (f)}{r^j},
$$
where
$\alpha = (\alpha _1 ,\alpha_2 ,...,\alpha _n ):\,\,\,\left| \alpha \right| = j \ge
2n,\,\,\,\alpha _p = 0$ if $r + ml_p = 0$ and $\,\alpha _p \ge 2$
if $k_p \ne 0,\,\,\,p = 1,2,...,n$. If $r + ml_p = \alpha _p = 0$,
we let $\left| {r + ml_p } \right|^{\alpha _p }=1$.

By the definition of $\tau _f (r)$ there is a number $r_0 (f)$ such that for
all $r > r_0 (f)$ we have
\[
\begin{array}{l}
 \tau _f (r) = \min \{M_0 (f),\,\,\frac{M_1 (f)}{r},\,\frac{M_2
(f)}{r^2},...,\frac{M_{2n} (f)}{r^{2n}}, \mathop {\inf }\limits_{j
\ge 2n + 1} \frac{M_j (f)}{r^j}\} = \mathop {\inf }\limits_{j \ge
2n + 1}
\frac{M_j (f)}{r^j}. \\
 \end{array}
\]
Hence there is a constant $C_1 (f)$ such that $\left| {a_{m,r} } \right| \le
C_1 (f)\tau _f (r)$ and

\begin{equation}
\label{eq9}
\left| {L_m (f,z)} \right| \le C_1 (f)(1 + \sum\limits_{r = 1}^{m - 1} {\tau
_f (r)t^{nr})}
\end{equation}

\noindent
for all $z = (z_1 ,z_2 ,....,z_n ):\,\,\,\,\,\,\frac{1}{t} \le \left| {z_p }
\right| \le t,\,\,\,\,p = 1,2,...,n$.

Now we estimate following expression on $T^n$:

\[
\frac{f(z^0) - L(f;z^0)}{(z_1^0 )^m + (z_2^0 )^m + ... + (z_n^0 )^m - n}.
\]
It is easy to show that
\[
\left| {\frac{(z_1^0 )^{m\beta _1 l_1 }(z_2^0 )^{m\beta _2 l_2 }....(z_n^0
)^{m\beta _n l_n } - 1}{(z_1^0 )^m + (z_2^0 )^m + ... + (z_n^0 )^m - n}}
\right| \le \left| l \right|,\,\,\,\,\,\,\forall z^0 \in T^n,\,\,\beta \in
\Lambda ,\,\,\,\forall l \in {\rm Z}_ + ^n .
\]
Moreover, the expansion (\ref{eq1}) can be rewritten as follows:
\[
f(z^0) = \sum\limits_{k \in {\rm Z}}^n(c_k (z^0))^k = \sum\limits_{l \in
{\rm Z}_ + ^n } \sum\limits_{r = 0}^{m - 1} \sum\limits_{\beta \in \Lambda
} \big(c_{\beta r + m\beta e} (z^0)\big)^{\beta r + m\beta e}.
\]
From here, using (\ref{eq7}) at $\left| \alpha \right| = j \ge 3n$, we obtain
\[
\left| {\frac{f(z^0) - L(f;z^0)}{(z_1^0 )^m + (z_2^0 )^m + ... + (z_n^0 )^m
- n}} \right| \le
\]

\[
 \le \left| {\sum\limits_{l \in {\rm Z}_ + ^n } {\sum\limits_{r = 0}^{m - 1}
{\sum\limits_{\beta \in \Lambda } {c_{\beta r + m\beta l} \frac{(z_1^0
)^{\beta _1 r + m\beta _1 l}(z_2^0 )^{\beta _2 r + m\beta _2 l}...(z_n^0
)^{\beta _n r + m\beta _n l}}{(z_1^0 )^m + (z_2^0 )^m + ... + (z_n^0 )^m -
n}} } } } \right| \le
\]

\[
 \le \sum\limits_{l \in {\rm Z}_ + ^n } {\sum\limits_{r = 0}^{m - 1}
{\sum\limits_{\beta \in \Lambda } {\left| l \right|\left| {c_{\beta r +
m\beta e} } \right| \le 2^nM_j (f)\sum\limits_{l \in {\rm Z}_ + ^n }
{\sum\limits_{r = 0}^{m - 1} {\frac{\left| l \right|}{\left| {r + ml_1 }
\right|^{\alpha _1 }\left| {r + ml_2 } \right|^{\alpha _2 }...\left| {r +
ml_n } \right|^{\alpha _n }}} } } } } .
\]
Consequently,

\[
\begin{array}{l}
 \left| {\frac{f(z^0) - L_m (f;z^0)}{(z_1^0 )^m + (z_2^0 )^m + ... + (z_n^0
)^m - n}} \right| \le 2^n \cdot m \cdot M_j (f)\sum\limits_{l \in
{\rm Z}_ + ^n } {\frac{\left| l \right|}{m^jl^\alpha } \le } 4^n \cdot n \cdot m\frac{M_j (f)}{m^j},\,\,j \ge 3n \\
 \end{array},
\]

\noindent
where $\alpha = (\alpha _1 ,\alpha _2 ,...,\alpha _n ):\,\,\,\left| \alpha
\right| = j \ge 3n,\,\,\,\alpha _p = 0$ if $r + ml_p = 0$ and $\,\alpha _p
\ge 3$ if $k_p \ne 0,\,\,\,p = 1,2,...,n$. If $r + ml_p = \alpha _p = 0$, we
let

\[
\left| {r + ml_p } \right|^{\alpha _p } = 1.
\]
By the definition of $\tau _f (r)$ there is a number $m_0 = m_0 (f)$ such
that for any
$m > m_0$ we have $\tau _f (m) = \mathop {\inf }\limits_{j \ge 3n + 1}
\frac{M_j (f)}{m^j}$.
Consequently, at $z:\,\,\,\left| {z_p } \right| < t,\,\,(t > 1),\,\,\,p =
1,2,...,n\,$ we obtain

\[
\begin{array}{l}
 \left| {\frac{z_1^m + z_2^m + ... + z_n^m - n}{(z_1^0 )^m + (z_2^0 )^m +
... + (z_n^0 )^m - n}(f(z^0) - L_m (f,z^0))} \right| \le 2 \cdot
4^n \cdot n^2 \cdot
m\tau _f (m)t^m,\,\,\quad m > m_0. \\
 \end{array}
\]
Hence there is some constant $C_2 (f)$ such that

\begin{equation}
\label{eq10} \left| {\frac{z_1^m + z_2^m + ... + z_n^m - n}{(z_1^0
)^m + (z_2^0 )^m + ... + (z_n^0 )^m - n}(f(z^0) - L_m (f,z^0))}
\right| \le C_2 (f)(1 + m\tau _f (m)t^m),\,\,\,\,\,\]
 \[ \forall z^0 \in T^n,\,\,\forall
z:\,\,\,\left| {z_p } \right| < t,\,\,p = \overline {1,n}
,\,\,\,\,\forall m \ge 1
\end{equation}

From inequalities (\ref{eq9}) and (\ref{eq10}) we obtain
inequality (\ref{eq8}). Using formula (\ref{eq4}) we define
following sequence $t_m $:
\[
t_m = \min \{\frac{1}{(r^3\tau _f (r))^{1 / nr}}:\,\,\,\,1 \le r \le m\}.
\]
Clearly, the sequence $t_m $ is decreasing and according to
(\ref{eq3}) we also have $t_m > 1$ for all $m$. Moreover,
\\$r^3\tau _f (r)t_m^{nr} \le 1$ for all $r \le m$.
If $z = (z_1 ,z_2 ,...,z_n ):\,\,\,\,\,\textstyle{1 \over {t_m }} \le \left|
{z_p } \right| \le t_m ,\,\,\,\,p = 1,2,...,n$, then inequality (\ref{eq9}) implies
that
\[
\left| {L_m (f;z^0;z)} \right| \le C_f (1 + \sum\limits_{r = 1}^m
{\frac{r^3\tau _f (r)t_m^{nr} }{r^2}} ) \le C_f (1 + \sum\limits_{r = 1}^m
{\frac{1}{r^2}} ) \le 3C_f .
\]
Thus, we have defined a mapping
\[
F_m (z) = (z,L_m (f;z^0;z))
\]
in a domain
\[
D(t_m ) = \{z \in {\rm C}^n:\,\,\,\textstyle{1 \over {t_m }} \le \left| {z_p
} \right| \le t_m ,\,\,\,p = 1,2,...,n\},
\]
whose image (along with graph of function $f)$ belongs to some ball $B(0,R_f
)$ from ${\rm C}^{n + 1}$ centered at the origin and with radius $R_f $ for any
$m \in {\rm N},\,\,\,\,z^0 \in T^n$. Moreover it has following properties:

1) $F_m (z^{(l)}) = (z^{(l)},f(z^{(l)})) \in \Gamma _f ,\,\,\,z^{(l)} =
(e^{\textstyle{{2\pi il_1 } \over m}},e^{\textstyle{{2\pi il_2 } \over
m}},...,e^{\textstyle{{2\pi il_n } \over m}}),\,\,\,l \in J_m $;

2) $F_m (z^0) = (z^0,f(z^0)) \in \Gamma _f ,\,\,\,z^0 \in T^n$.

According to the hypothesis of the proposition we have
\[
\overline {\mathop {\lim }\limits_{m \to \infty } } m^{\textstyle{1 \over {n
+ 1}}}\ln t_m = + \infty .
\]
From the sequence $m^{\textstyle{1 \over {n + 1}}}\ln t_m $ we can choose a
subsequence which tends to infinity$. $Without loss of generality we can assume
that
\[
\mathop {\lim }\limits_{m \to \infty } m^{\textstyle{1 \over {n + 1}}}\ln
t_m = + \infty .
\]
Next, we consider weighted multipole Green's functions $G_m (w)$
in a ball $B(0,2R_f )$ with poles $w^{(l)} =
(z^{(l)},f(z^{(l)}))\,,\,\, \\l \in J_m $ and with weight $m^{ -
\textstyle{n \over {n + 1}}}$ at each pole (see [8], [9]):
\[
\begin{array}{l}
 G_m (w) = \sup \left\{ {u(w) \in psh(B(0,2R_f )):\,\,u(w) \le 0,} \right.
\\
 \left. {u(w) - m^{ - \textstyle{n \over {n + 1}}}\ln \left| {w - w^{(l)}}
\right| = O(\ref{eq1}),\,\,w \to w^{(l)},\,l \in J_m } \right\} \\
 \end{array}
\]
and $(dd^cG_m (w))^{n + 1} = 0$ on $B(0,2R_f )\backslash \left\{
{w^{(l)},\,\,\,\,l \in J_m } \right\}$, $G_m (w) = 0$ on $\partial B(0,2R_f
)$, $G_m (w) - m^{ - \textstyle{n \over {n + 1}}}\ln \left| {w - w^{(l)}}
\right| = O(\ref{eq1}),\,\,w \to w^{(l)},\,l \in J_m $. Moreover,
\[
(dd^cG_m )^{n + 1} = (2\pi )^{n + 1}\sum\limits_{l \in J_m } {m^{ - n}\delta
_{w^{(l)}} } ,
\]
where $\delta _{w^{(l)}} $ is the Dirac mass at $w^{(l)}$. From here we
obtain that
\[
\int\limits_{B(0,2R_f )} {(dd^cG_m (w))} ^{n + 1} = (2\pi )^{n + 1}.
\]
We consider following plurisubharmonic function in a domain $D(t_m )$:
\[
u_m (z) = G_m (F_m (z)),
\]
which at a point $z^0 \in T^n$ takes on the value $G_m (w^0) = G_m
(z^0,\,f(z^0))$ and at a point $z^{(l)},\,\,l \in J_m $, has poles with
weights $m^{ - \textstyle{n \over {n + 1}}}$.

Let $g(\lambda )$ be the Green's function on a ring
\[
U_m = \left\{ {\lambda \in {\rm C}:\,\,\,\frac{1}{t_m^m } < \left| \lambda
\right| < t_m^m } \right\}
\]
which at a point $\lambda = 1$ has a pole with weight $1$, such that
\[
\mathop {\sup }\limits_{\left| \lambda \right| = 1} g(\lambda ) \le C_g \log
t_m^m ,
\]
where $C_g $ is some negative constant (existence of such a Green function
is proved in [5]).

Let
\[
V(z_1 ,z_2 ,...,z_n ) = \mathop {\max }\limits_{1 \le j \le n} g(z_j
),\,\,\,\,(z_1 ,z_2 ,...,z_n ) \in D(t_m^m ).
\]
Then $V(z_1 ,z_2 ,...,z_n )$ is a maximal plurisubharmonic function with a
single pole at a point $a = (1,1,...,1) \in T^n$. The function $V(z_1^m
,z_2^m ,...,z_n^m )$ is also maximal plurisubharmonic function with poles at
points $z^{(l)},\,\,l \in J_m $. Then for all $z \in D(t_m )$ we have the
inequality
\[
u_m (z) \le m^{ - \textstyle{n \over {n + 1}}}V(z_1^m ,z_2^m ,...,z_n^m ).
\]
We note that $t_m^m \to + \infty $ and, consequently, $t_m^m > 2$ for
sufficiently large $m$. From here we obtain
\[
G_m (z^0,\,f(z^0)) = u_m (z^0) \le C \cdot m^{ - \textstyle{n \over {n +
1}}}\ln t_m^m = C \cdot m^{1 - \textstyle{n \over {n + 1}}}\ln t_m ,
\]
where $C < 0$. Thus,
\[
\mathop {\lim }\limits_{m \to \infty } G_m (z^0,\,f(z^0)) = - \infty
,\,\,\,\,\,z^0 \in T^n.
\]
That is, $\Gamma _f $ is pluripolar, which completes the proof of the
proposition.

\vskip 0.5cm

\textbf{Proof of the theorem. }Let $f(z) \in C^\infty (T^n)$ be a Denjoy
quasianalytic function. We note that the graph of the function $f$
is pluripolar if and only if the graph of $cf$ is pluripolar, where $c \ne
0$ is a constant. Multiplying $f$ by a small constant we can assume that our
smooth function $f$ satisfies the inequality
\[
M_3 (f) < \frac{1}{2}.
\]
Let

\begin{equation}
\label{eq11}
\tilde {\tau }_f (r) = \mathop {\inf }\limits_{s \ge 3} \frac{M_s (f)}{r^{s
- 3}} = \mathop {\inf }\limits_{s \ge 0} \frac{M_{s + 3} (f)}{r^s} <
\frac{1}{2}
\end{equation}
be the associated function for the shifted sequence $\left\{ {\tilde {M}_s }
\right\} = \left\{ {M_{s + 3} (f)} \right\}$. Setting

\begin{equation}
\label{eq12}
\ln \theta _f (m) = \min \left\{ { - \frac{\ln \tilde {\tau }_f
(r)}{nr}:\,\,1 \le r \le m} \right\},
\end{equation}
we obtain

\begin{equation}
\label{eq13}
\ln t_m \ge \ln \theta _f (m) > 0.
\end{equation}
From the definition of $\tau _f (r)$ we have
\[
r^3\tau _f (r) = \min \left\{ {M_0 (f)r^3,\,\,M_1 (f)r^2,\,\,M_2
(f)r,\,\,\mathop {\inf }\limits_{s \ge 3} \frac{M_s (f)}{r^{s - 3}}}
\right\} = \tilde {\tau }_f (r)
\]
for all $r > r_0 (f) \ge 0$. We will show that

\begin{equation}
\label{eq14}
\mathop {\overline {\lim } }\limits_{m \to \infty } m^{\textstyle{1 \over {n
+ 1}}}\ln \theta _f (m) = \infty .
\end{equation}

To show this we use the following lemma from [5].

\textbf{Lemma. } \textit{Let}  $\tilde {h}(s) = h(e^s)$\textit{ - positive, increasing, convex function of variable }$s$\textit{ on }$[0;\infty )$\textit{ and let }

\[
\tilde {H}(s) = \min \left\{ {\tilde {h}(s)e^{ - s}:\,\,\,0 \le s \le x}
\right\}.
\]

\textit{If }$\tilde {H}(x) \le Ce^{ - \alpha x},\,\,0 < \alpha < 1$,\textit{ for all }$x \ge 0$\textit{, then}

\[
\int\limits_0^\infty {\tilde {h}(s)e^{ - s}ds = \int\limits_1^\infty
{\frac{h(t)}{t^2}} } dt < \infty .
\]

Suppose that condition (\ref{eq14}) does not hold, that is,

\[
\mathop {\overline {\lim } }\limits_{m \to \infty } m^{\textstyle{1 \over {n
+ 1}}}\ln \theta _f (m) \le C < \infty .
\]
Let $h(t) = - \ln \tilde {\tau }_s (t)$ and $\tilde {h}(s) = h(e^s)$. We
have
\[
H(t) = \ln \theta _f (t) = \min \left\{ {\textstyle{{h(r)} \over
{nr}}:\,\,\,1 \le r \le t} \right\}.
\]
It follows that $H(m) < C / m^{\textstyle{1 \over {n + 1}}}$ for each
positive $m$. Since $H$ is decreasing, we have $H(x) \le 2C / x^{\textstyle{1
\over {n + 1}}}$ for all $x$.
Computations yield
\[
\tilde {H}(x) = \min \left\{ {\tilde {h}(s)e^{ - s}:\,\,0 \le s \le x}
\right\} = H(e^x),\,\,\,H(x) < 2Ce^{ - \textstyle{x \over {n +
1}}},\,\,\forall x \ge 0.
\]
According to (\ref{eq12}) we have $h(t) = - \ln \tilde {\tau }_s (t) > 0$ and using
the lemma we have
\[
\int\limits_1^\infty {\frac{h(t)}{t^2}dt < \infty } .
\]
On the other hand, we have $r^3\tau _f (r) = \tilde {\tau }_f (r)$ for all
$r > r_0 (f) \ge 0$ and the function $f$ is a Denjoy quasianalytic function.
Then according to Carleman's theorem we have
\[
\int\limits_1^\infty {\frac{ - \ln \tilde {\tau }_f (r)}{r^2}dr = }
\int\limits_1^\infty {\frac{ - 3\ln r}{r^2}dr + \int\limits_1^\infty {\frac{
- \ln \tau _f (r)}{r^2}dr = } } + \infty ,
\]
which is a required contradiction. Therefore, the function $f$ satisfies
condition (\ref{eq14}) and, consequently, it also satisfies condition (\ref{eq3}). Thus, any
Denjoy quasianalytic function satisfies the conditions of the proposition.
It follows that $\Gamma _f $ is pluripolar. The proof of the theorem is
complete.

\begin{center}
\textbf{REFERENCES }
\end{center}

\begin{enumerate}
\item Diederich K., Fornass J.E. A smooth curve in ${\rm {\bf
C}}^{\rm {\bf 2}}$ which is not a pluripolar set. Duke Math. J.
1982, Volume 49, 931-936.

\item Diederich K., Fornass J.E. Smooth, but not complex-analytic
pluripolar sets. Manuscript Math. 1982, Volume 37, 121-125.

\item Edigaryan A., Weegerinck J. Graphs that are not complete pluripolar.
Proceedings of the AMS. Volume 131, ¹ 8, 2459-2465

\item Edlund T., Joricke B. The pluripolar hull of a graph and fine analytic
continuation. Ark. Mat., Volume\textbf{ 44 }(2006), 39--60.

\item Coman D., Levenberg N., Poletskiy E.A. Quasianalyticity and Pluripolarity.
J. Amer. Math. Soc. Volume 18, (2005) ¹2, 9-16.

\item Shcherbina N. Pluripolar graphs are holomorphic// Acta Math., Volume 194
(2005), 203-216.

\item Mandelbrojt S. Quasianalytic class of functions. Moscow, 1937.

\item Demailly J.P. Measures de Monge-Ampere et mesures plurisousharmoniques//
Math. Z. Volume\textbf{ 194,} (1987), 519-564.

\item Lelong P. Function de Green pluricomplexe et lemmes de Schwarz dans les
espaces de Banach// J. Math. Pures Appl. Volume\textbf{ 68,}
(1989), 319-347.
\end{enumerate}

Sevdiyor A. Imomkulov

Department of Mathematics and Physics,\\ Navoi State Pedagogical
Institute (Uzbekistan)\\e-mail: \underline {sevdi@rambler.ru}\\

Zafar Sh. Ibragimov \\ Department of Mathematics and
Physics,\\Urgench State University (Uzbekistan)\\e-mail:
\underline {z.ibragim@gmail.com}

\end{document}